\documentclass[12pt,leqno]{article}
\usepackage{amsopn,amsmath,amssymb,amscd,latexsym}
\usepackage[all]{xy}
\newcommand{\C}{{\mathbb{C}}}
\newcommand{\F}{{\mathbb{F}}}
\newcommand{\oF}{\overline{\F}}
\newcommand{\G}{\mathbb{G}}
\newcommand{\Q}{{\mathbb{Q}}}
\newcommand{\R}{{\mathbb{R}}}
\newcommand{\Z}{{\mathbb{Z}}}
\newcommand{\Co}{\mathrm{Co}}
\newcommand{\ddet}{\mathrm{det}}
\newcommand{\oL}{\overline{L}}
\newcommand{\loc}{\mathrm{loc}}
\newcommand{\oM}{\overline{M}}
\renewcommand{\mod}{\;\mathrm{mod}\;}
\newcommand{\oR}{\overline{R}}
\newcommand{\red}{\mathrm{red}}
\newcommand{\sgn}{\mathrm{sgn}\,}
\newcommand{\spec}{\mathrm{spec}\,}
\newcommand{\tors}{\mathrm{tors}}
\newcommand{\Arg}{\mathrm{Arg}\,}
\newcommand{\card}{\mathrm{card}\,}

\newcommand{\End}{\mathrm{End}\,}
\newcommand{\Harm}{\mathrm{Harm}}

\newcommand{\Ind}{\mathrm{Ind}}
\newcommand{\Ker}{\mathrm{Ker}\,}
\newcommand{\Log}{\mathrm{Log}}
\newcommand{\RRe}{\mathrm{Re}\,}
\newcommand{\Tr}{\mathrm{Tr}}
\newcommand{\vol}{\mathrm{vol}}
\newcommand{\Ah}{{\mathcal A}}

\newcommand{\Dh}{{\mathcal D}}
\newcommand{\Eh}{{\mathcal E}}
\newcommand{\Fh}{{\mathcal F}}
\newcommand{\FL}{\Fh\Lh}
\newcommand{\Lh}{{\mathcal L}}
\newcommand{\Oh}{{\mathcal O}}

\newcommand{\emm}{{\mathfrak{m}}}
\newcommand{\oemm}{\overline{\emm}}

\newcommand{\oxi}{\overline{\xi}}
\newcommand{\opi}{\overline{\pi}}
\newcommand{\opartial}{\overline{\partial}}
\newcommand{\ohne}{\setminus}
\newcommand{\silo}{\stackrel{\sim}{\longrightarrow}}
\newcommand{\tei}{\, | \,}
\newcommand{\hullet}{\raisebox{0.05cm}{$\scriptscriptstyle \bullet$}}
\newcommand{\halb}{\frac{1}{2}}
\newcommand{\dis}{\displaystyle}

\newtheorem{theorem}{Theorem}[section]
\newtheorem{prop}[theorem]{Proposition}
\newtheorem{cor}[theorem]{Corollary}

\newenvironment{rems}{\noindent {\bf Remarks}}{}
\newenvironment{proofof}{\noindent {\bf Proof of}}{\mbox{}\hfill$\Box$}
\newenvironment{proof}{\noindent {\bf Proof}}{\mbox{}\hfill$\Box$}
\begin{document}
\title{On the nature of the ``explicit formulas'' in analytic number theory --- a simple example} 
\author{Christopher Deninger}
\date{\ }
\maketitle

\section{Introduction}
\label{sec:1}

The ``explicit formulas'' in analytic number theory relate the prime numbers with the zeroes of the Riemann zeta function. They were first used by Riemann in his famous note which started analytic number theory. Various authors, notably A. Weil, subsequently developed generalizations of these formulas. A simple version of the ``explicit formulas§ for the Riemann zeta function is the following. Given a test function $\alpha \in C^{\infty}_0 (\R)$ we have
\begin{eqnarray}
  \label{eq:1}
  \Phi (0) - \sum_{\rho} \Phi (\rho) + \Phi (1) & = & \sum_p \log p \sum_{k \ge 1} \alpha (k \log p)\\ 
& &  + \sum_p \log p \sum_{k \leq -1} p^k \alpha (k \log p) + W_{\infty} (\alpha) \; . \nonumber
\end{eqnarray}
Here $\rho$ runs over the non-trivial zeroes of $\zeta (s)$ and $\Phi$ is the entire function defined by the integral
\[
\Phi (s) = \int_{\R} e^{ts} \alpha (t) \, dt \; .
\]
Moreover we have set
\[
W_{\infty} (\alpha) = \alpha (0) \log \pi + \int^{\infty}_0 (\alpha (t) + e^{-t} \alpha (-t)) (1 - e^{-2t})^{-1} - \alpha (0) t^{-1} e^{-2t} \, dt \; .
\]
Formula (1) and generalizations can be found in \cite{B} for example. We will view (1) as the ``explicit formula'' for the function
\[
\hat{\zeta} (s) = \zeta (s) \zeta_{\infty} (s) \; , \; \mbox{where} \quad \zeta_{\infty} (s) = \pi^{-s/2} \Gamma \left( \frac{s}{2} \right) \; .
\]
In a sense $\hat{\zeta} (s)$ is an Euler product over {\it all} valuations of $\Q$. Note that $\hat{\zeta} (s)$ has simple poles at $s = 0,1$ and that its zeroes are exactly the non-trivial zeroes of $\zeta (s)$. 

There are good reasons to think that formula (1) has conceptual origins rooted in analysis on certain foliated spaces. These spaces however remain to be found e.g. \cite{D1} , \cite{D2}. 

In this note we will be concerned with ``explicit formulas'' for certain function fields.\\
A function field $F$ (of transcendence degree one) over  the finite field $\F_q$ with $q = p^r$ elements is a finite extension of the field $\F_q (T)$ of rational functions over $\F_q$ such that $\F_q$ is algebraically closed in $F$. Since the time of E. Artin's thesis it is well known that there are many analogies between the arithmetic of number fields and such function fields. In particular a zeta function can be defined for $F$, 
\[
\hat{\zeta}_F (s) = \prod_w (1 - Nw^{-s})^{-1} \; .
\]
Here $w$ runs over the valuations of $F$ and $Nw$ is the number of elements in the (finite) residue field of $w$. This Euler product converges for $\RRe s > 1$ and $\hat{\zeta}_F (s)$ can be meromorphically continued to the entire plane. In fact, it is a rational function of $q^{-s}$. The function $\hat{\zeta}_F (s)$ has simple poles at the following points on $\RRe s = 0$ resp. $\RRe s = 1$:
\[
\frac{2\pi i\nu}{\log q} \quad \mbox{resp.} \quad 1 + \frac{2\pi i\nu}{\log q} \quad \mbox{for} \; \nu \in \Z \; .
\]
Using geometric methods A. Weil was able to prove that the Riemann hypotheses holds for $\hat{\zeta}_F (s)$ i.e. that all zeroes $\rho$ of $\hat{\zeta}_F (s)$ lie on the line $\RRe s = 1/2$. 

The usual methods give the following analogue of (1) if $F$ has genus $g$:
\begin{eqnarray}
  \label{eq:2}
 \lefteqn{\hspace*{0.5cm} \sum_{\nu \in \Z} \Phi \left( \frac{2 \pi i \nu}{\log q} \right) - \sum_{\rho} \Phi (\rho) + \sum_{\nu \in \Z} \Phi \left( 1 + \frac{2 \pi i \nu}{\log q} \right) } \\
& =& \alpha (0) (2 - 2g) \log q + \sum_w \log Nw \sum_{k \ge 1} \alpha (k \log Nw) \nonumber \\
 & & + \sum_w \log Nw \sum_{k \le -1} Nw^k \alpha (k \log Nw) \; . \nonumber
\end{eqnarray}
Note that there is no term like $W_{\infty} (\alpha)$ which corresponds to the Archimedian valuation of $\Q$. In the function field case all valuations are non-Archimedian.

Now let $E_0$ be an elliptic curve over $\F_q$ with function field $F = \F_q (E_0)$. In this note we will interpret formula (2) for $F = \F_q (E_0)$ as a transversal index theorem for an $\R$-action on a suitable $3$-dimensional space $X$ which is not a manifold. This is the content of section 3.

Related spaces for elliptic curves have been considered by Ihara but not published. See \cite{I1}, \cite{I2}, \cite{I3}, \cite{I4} however for his constructions in the cyclotomic, modular and Shimura curve cases.

 In section 2 we provide some background on classical transversal index theory on manifolds. Finally in the last section we give another example of a transversal index calculation on our space $X$ above.

In the field of arithmetic topology, see e.g. \cite{Sik}, certain analogies are studied between number fields and $3$-manifolds. The constructions of the present note make some of these analogies concrete in the easy case where the number field is replaced by $F = \F_q (E_0)$. See for instance the explicit bijection in proposition 3.3 between valuations of $F$ and certain knots in the $3$-space $X$.

I would like to thank Y. Ihara, Y. Kordyukov, P. Schneider and M. Volkov for helpful comments and discussions. I am grateful to S. Kanemitsu, Professor of mathematics and head of the international Kanemitsu corporation \cite{R} for the invitation to the joint Japanese--Chinese meeting in analytic number theory. I would also like to thank K. Matsumoto for kindly supporting my stay in Japan. Finally I am very grateful to the referee for valuable suggestions.

\section{A short introduction to transversal Index theory}
\label{sec:2}

Transversal index theory is concerned with (pseudo-)differential operators or complexes of such operators, which are elliptic in the directions transversal to the orbits of an action by a Lie group $G$.

For compact groups the theory was initiated in Atiyah's lecture note \cite{A}. For non-compact groups the definition of the transversal index as a distribution on $G$ is due to H\"ormander \cite{Si}, \cite{NZ}. Connes and Moscovici have greatly generalized the theory using non-commutative methods.

Let us recall H\"ormander's definition of the transversal index. Let $X$ be a compact manifold and consider a complex of differential operators acting on the smooth sections of vector-bundles $E_0 , \ldots , E_N$ over $X$:
\begin{equation}
  \label{eq:3}
  0 \longrightarrow C^{\infty} (E_0) \xrightarrow{D_0} C^{\infty} (E_1) \xrightarrow{D_1} \ldots \longrightarrow C^{\infty} (E_N) \longrightarrow 0 \; .
\end{equation}
Let us assume that a Lie group $G$ acts smoothly on $X$ and that we are given smooth (co-)actions $\tilde{g}_j : g^* E_j \to E_j$ on the bundles $E_j$. Furthermore the differential operators should commute with the induced $G$-actions on the $C^{\infty} (E_j)$. Let $\pi : T^* X \ohne 0 \to X$ denote the projection. Then associated to (3) there is the symbol sequence -- a complex of vector bundles over $T^* X \ohne 0$:
\begin{equation}
  \label{eq:4}
  0 \longrightarrow \pi^* E_0 \xrightarrow{\sigma (D_0)} \pi^* E_1 \longrightarrow \ldots \longrightarrow \pi^* E_N \longrightarrow 0 \; .
\end{equation}
Define the characteristic variety at the $j$-th place $C_j \subset T^* X \ohne 0$ to be the set of $\xi$ in $T^* X \ohne 0$ over which (4) is not exact.

The complex (3) is called transversally elliptic with respect to the $G$-action, if all $C_j$ are disjoint with $T^*_G X \ohne 0$. Here $T^*_G X$ is the closed $G$-invariant subspace of $T^*X$ consisting of those $\xi \in T^* X$ that annihilate all vector fields along the $G$-orbits \cite{A} p. 7.

Choose hermitian metrics on $E_j$ and a volume form on $X$ (or density if $X$ is non-orientable). This gives a scalar product on each $C^{\infty} (E_j)$. Let $L^2 (E_j)$ be the $L^2$-completion. We extend the operator
\[
D_j : C^{\infty} (E_j) \longrightarrow C^{\infty} (E_{j+1})
\]
to a closed densely defined operator
\[
\tilde{D}_j : L^2 (E_j) \longrightarrow L^2 (E_{j+1})
\]
by setting $\tilde{D}_j = (D^{\dagger}_j)^*$ where $D^{\dagger}_j$ is the formal adjoint of $D_j$. The domain of definition of $\tilde{D}_j$ consists of those $f \in L^2 (E_j)$ for which ${D}_j f$, taken in the distributional sense is in $L^2 (E_{j+1})$. The operator $\tilde{D}_j$ is known as the weak (or maximal) closure of $D_j$.

For simplicity assume that $G$ acts conformally on the $L^2 (E_j)$. Then $G$ commutes with $D^{\dagger}_j , \tilde{D}_j$ and hence it acts on
\[
\Harm^j_{L^2} (X) = \Ker \tilde{D}_j \cap \Ker (\tilde{D}_j)^* \; .
\]
As a closed subspace of $L^2 (E_j)$ this is a Hilbert space. For $\alpha \in C^{\infty}_0 (G)$ define an endomorphism of $\Harm^j_{L^2} (X)$ by the formula:
\[
S_j (\alpha) = \int_G \alpha (g) g^* \, dg \; .
\]
Then H\"ormander proved, c.f. \cite{Si} p. 29, \cite{NZ}:

\begin{theorem}
  \label{t21}
$S_j (\alpha)$ is of trace class and the map $\alpha \mapsto \Tr (S_j (\alpha))$ is a distribution $T_j$ on $G$. The representation of $G$ on $\Harm^j_{L^2} (X)$ decomposes into a discrete direct sum of irreducible representations. 
\end{theorem}
The transversal index of (3) can now be defined as the distribution:
\[
\Ind_t (D) = \sum^N_{j=0} (-1)^j T_j \quad \mbox{in} \; \Dh' (\R) \; .
\]
We go on to describe two computations of a transversal index for a $G = \R$-action due to \'Alvarez L\'opez, Kordyukov \cite{AK} and Lazarov \cite{L}.

Consider a compact $(d+1)$-dimensional manifold $X$ with a smooth one-codimensional foliation $\Fh$ and a smooth $\R$-action whose orbits are everywhere transversal to the leaves of $\Fh$ and which maps leaves to leaves.

Then the de Rham complex along the leaves of $\Fh$:
\[
0 \longrightarrow C^{\infty} (\Lambda^0 T^* \Fh) \xrightarrow{d_{\Fh}} C^{\infty} (\Lambda^1 T^* \Fh) \xrightarrow{d_{\Fh}} \ldots \longrightarrow C^{\infty} (\Lambda^d T^* \Fh) \longrightarrow 0
\]
is transversally elliptic to the $\R$-action. Under a certain technical condition -- $\R$ has to act isometrically on $T\Fh$ -- the transversal index $\Ind_t (d_{\Fh})$ is defined as above.

We assume that all compact orbits $\gamma$ of the $\R$-action are non-degenerate in the sense that for all $x \in \gamma$ and $k \in \Z \ohne 0$ the endomorphism $T_x \phi^{kl (\gamma)}$ has the eigenvalue $1$ with algebraic multiplicity one. Here $l (\gamma)$ is the length of $\gamma$. Note that the tangent vector $Y_x$ to the $\R$-orbit through $x$ is kept fixed by $T_x \phi^{kl (\gamma)}$. 

We set
\[
\varepsilon_{\gamma} (k) = \sgn \det (1 - T_x \phi^{kl (\gamma)} \tei T_x X / \R \cdot Y_x) = \sgn \det (1 - T_x \phi^{kl (\gamma)} \tei T_x \Fh) \; .
\]
The following theorem is proved in \cite{AK}, \cite{L} and (on $\R^*$) in \cite{DS}:

\begin{theorem}
  \label{t22}
$\dis \Ind_t (d_{\Fh}) = \chi_{\Co} (\Fh , \mu) \delta_0 + \sum_{\gamma} l (\gamma) \sum_{k \in \Z \ohne 0} \varepsilon_{\gamma} (k) \delta_{kl (\gamma)}$ in $\Dh' (\R)$.
\end{theorem}
Here $\chi_{\Co} (\Fh , \mu)$ is Connes' Euler characteristic of the foliation with respect to a certain transversal measure $\mu$.

\begin{rems}
 1) The version of this theorem in \cite{AK} uses a more general definition of the transversal index and does not require that $\R$ acts isometrically on $T\Fh$. \\
2) We will not discuss the term $\chi_{\Co} (\Fh , \mu)$ because in the situation of the next section its analogue vanishes.
\end{rems}

Let us now assume that $d = 2g$ is even and that $\Fh$ is actually a foliation by complex manifolds i.e. that there exists a smooth almost complex structure $J$ on $T\Fh$ whose restriction to every leaf is integrable. Moreover every $t \in \R$ should induce (bi-)holomorphic maps between leaves. Then we can consider the Dolbeault complex along the leaves of $\Fh$:
\[
0 \longrightarrow C^{\infty} (\Lambda^0_{\C} T^*_c \Fh) \xrightarrow{\overline{\partial}_{\Fh}} C^{\infty} (\Lambda^1_{\C} T^*_c \Fh) \xrightarrow{\overline{\partial}_{\Fh}} \ldots \longrightarrow C^{\infty} (\Lambda^g_{\C} T^*_c \Fh) \longrightarrow 0 \; .
\]
The following result is a special case of Lazarov's theorem 2.10 in \cite{L}:

\begin{theorem}
  \label{t23}
In $\Dh' (\R)$ the following equality holds:
\[
\Ind_t (\overline{\partial}_{\Fh}) = \chi_{\Co} (\Fh , \Oh, \mu) \cdot \delta_0 + \sum_{\gamma} l (\gamma) \sum_{k \in \Z \ohne 0} \ddet_{\C} (1 - T_{cx} \phi^{kl (\gamma)} \tei T_{cx} \Fh)^{-1} \delta_{kl (\gamma)} \; .
\]
Here $x$ is any point on the relevant orbit $\gamma$ and $\chi_{\Co} (\Fh , \Oh , \mu)$ is the holomorphic Connes' Euler characteristic of the foliation with respect to $\mu$.
\end{theorem}

\section{A transversal index calculation on a generalized solenoid}
\label{sec:3}
Consider an elliptic curve $E_0$ over $\F_q$, where $q = p^r$, together with its $q$-th power Frobenius endomorphism $\varphi_0 : E_0 \to E_0$ over $\F_q$.

We consider a lift of the pair $(E_0 , \varphi_0)$ to characteristic zero as follows.

First assume that $E_0$ is ordinary and let $R = W (\F_q)$ denote the ring of Witt vectors of $\F_q$ with maximal ideal $\emm = p R$ and quotient field $L$. Note that $R /\emm = \F_q$. Up to isomorphism there is then a unique elliptic curve $\Eh$ over $\spec R$ such that $\Eh \otimes \F_q = E_0$ and such that
\[
\End_R (\Eh) \silo \End_{\F_q} (E_0)
\]
is an isomorphism. This $\Eh$ is called the canonical lift of $E_0$.

In particular, there is a unique endomorphism $\varphi : \Eh \to \Eh$ lifting $\varphi_0$ i.e. such that
\[
(\Eh , \varphi) \otimes \F_q = (E_0 , \varphi_0) \; .
\]
Now let $E_0$ be supersingular. There exist the following:
\begin{itemize}
\item a complete local integral domain $R$ with field of fractions $L$ a finite extension of $\Q_p$ such that $R / \emm = \F_q$ if $\emm$ is the maximal ideal of $R$.
\item an elliptic curve $\Eh$ over $\spec R$ together with an endomorphism \\
$\varphi : \Eh \to \Eh$ such that
\[
(\Eh , \varphi) \otimes \F_q = (E_0 , \varphi_0) \; .
\]
\end{itemize}
The pair $(\Eh , \varphi)$ is not canonically determined. In the following, any choice will do if $E_0$ is supersingular. Proofs or references for these facts can be found in \cite{O}. Up to a finite extension of $\F_q$ they follow from \cite{Deu} p. 259--263.

Let $E_0 / \F_q$ be arbitrary and consider $(\Eh , \varphi)$ as above. We denote by $E = \Eh \otimes_R L$ the generic fibre. Then
\[
\End^0_L (E) = \End_L (E) \otimes \Q
\]
is a field $K$ which is either $\Q$ or an imaginary quadratic extension of $\Q$.

Fixing an embedding $L \subset \C$ we can consider the complex analytic elliptic curve $E (\C)$. Let $\omega$ be a non-zero holomorphic one-form on $E (\C)$ and let $\Gamma \subset \C$ be its period lattice. Then the Abel--Jacobi map
\[
E (\C) \silo \C / \Gamma \; , \; p \longmapsto \int^p_0 \omega \mod \Gamma
\]
defines an isomorphism.

There is a homomorphism:
\[
\Theta : K \cong \End^0_L (E) \longrightarrow \End^0 (\C / \Gamma) \; .
\]
Take the embedding $K \subset \C$ for which $\Theta (\alpha)$ induces multiplication by $\alpha$ on the Lie algebra $\C$ of $\C / \Gamma$. 

Let
\[
\xi \in \Theta^{-1} (\End_L (E)) \subset K \subset \C
\]
be the unique element with $\Theta (\xi) = \varphi \otimes L$. According to a theorem of Hasse -- the Riemann conjecture for elliptic curves over function fields -- we have $|\xi|^2 = q$.

It is clear that $\xi \Gamma \subset \Gamma$. Hence we may consider the ${\xi}$-adic ``Tate modules'':
\[
T_{\xi} \Gamma = \varprojlim \Gamma / \xi^n \Gamma \quad \mbox{and} \quad V_{\xi} \Gamma = T_{\xi} \Gamma \otimes \Q \; .
\]
Note that $V_{\xi} \Gamma$ is a $\Q_p$-vector space of dimension one or two.

We define an additive subgroup $V \subset \C$ by
\[
V = \Gamma [\xi^{-1}] = \bigcup_{\nu \ge 0} \xi^{-\nu} \Gamma \; .
\]
Then $V$ acts by translation on $V_{\xi} \Gamma$ since multiplication by $\xi$ is invertible on $V_{\xi} \Gamma$. We now introduce actions of $V$ and $\Lambda = (\log q) \Z \subset \R$ on the space
\[
\tilde{X} = \C \times V_{\xi} \Gamma \times \R \; .
\]
They are defined by the formulas:
\[
(z , \hat{v} , t) \cdot v = (z + v , \hat{v}- v, t) \quad \mbox{for} \; v \in V
\]
and
\[
(z , \hat{v} , t) \cdot \lambda = (\xi^{-\nu} z , \xi^{-\nu} \hat{v} , t + \lambda) \quad \mbox{for} \; \lambda = \nu \log q \in \Lambda \; .
\]
Together we get a free right action on $\tilde{X}$ of the semidirect product
\[
H = V \rtimes \Lambda = V \cdot \Lambda
\]
where $\lambda \cdot v = v^{\lambda} \cdot \lambda$ with $v^{\lambda} = \xi^{\nu} v$. 

The $3$-dimensional space on which we want to do transversal index theory is the compact quotient $X = \tilde{X} / H$. 

The Lie group in question is $\R$ which acts on $X$ by translation. For $s \in \R$ we set $\phi^s (z , \hat{v} , t) = (z , \hat{v} , t+s)$. In order to define the relevant differential operators we make the following observations. The space $X$ is a smooth foliated space in the sense of \cite{MS} in two ways. First there is the foliation $\Lh$ of $X$ by its $3$-dimensional path components. They are non-compact and $\R$-equivariantly diffeomorphic to either $\C \times \R$ or a suitable circle bundle over $\C$, see the end of this section for details.

The transversals for $\Lh$ are totally disconnected. Thus $X$ becomes a generalized solenoid in the sense of \cite{Su}. 

Secondly there is the foliation $\Fh\Lh$ of $X$ by the surfaces $\tau (\C \times \{ \hat{v}\} \times \{ t \}) $ where $\tau : \tilde{X} \to X$ is the natural projection. Its leaves are all diffeomorphic to $\C$. Any leaf $L$ of $\Lh$ is foliated as a $3$-dimensional manifold by the leaves $F$ of $\Fh \Lh$ that are contained in $L$. The resulting foliation $\Fh_L$ of $L$ is one-codimensional and everywhere transversal to the orbits of the $\R$-action restricted to $L$.

The union of the tangent spaces to the $\Fh \Lh$-leaves is naturally a rank $2$ vector bundle on $X$ denoted by $T \Fh \Lh$. As in \cite{MS} let $\Ah^i_{\Fh\Lh} (X)$ denote the space of those sections of $\Lambda^i T^* \Fh\Lh$ which are smooth on the $\Fh\Lh$-leaves and continuous transversally. We then have the de Rham complex along the $\Fh\Lh$-leaves \cite{MS}:
\begin{equation}
  \label{eq:5}
  0 \longrightarrow \Ah^0_{\Fh\Lh} (X) \xrightarrow{d_{\Fh\Lh}} \Ah^1_{\Fh\Lh} (X) \xrightarrow{d_{\Fh\Lh}} \Ah^2_{\Fh\Lh} (X) \longrightarrow 0 \; .
\end{equation}
Because of the equation
\[
(d_{\Fh\Lh} \omega) \, |_L = d_{\Fh_L} (\omega \, |_L) 
\]
this complex can be viewed as the family of de Rham complexes along the leaves $L$ of $\Lh$
\begin{equation}
  \label{eq:6}
  0 \longrightarrow \Ah^0_{\Fh_L}  (L) \xrightarrow{d_{\Fh_L}} \Ah^1_{\Fh_L} (L) \xrightarrow{d_{\Fh_L}} \Ah^2_{\Fh_L} (L) \longrightarrow 0 \; .
\end{equation}
These complexes are transversally elliptic with respect to the $\R$-action on $L$.

We now define the transversal index for the complex (5) similarly as in the case of manifolds.

Let $\mu_{\xi}$ denote a Haar measure on the locally compact abelian group $V_{\pi} \Gamma$. Then the measure 
\[
\tilde{\mu} = dx \, dy \otimes \mu_{\xi} \otimes dt
\]
on $\tilde{X}$ is $H$-invariant and hence induces a measure $\mu$ on $X$. In this regard, note that
\[
\mu_{\xi} (\xi^{\nu} A) = |\xi|^{-2\nu} \mu_{\xi} (A) = q^{-\nu} \mu_{\xi} (A) \; .
\]
It suffices to check this for $\nu \ge 0$ where it follows from the facts that $(\xi^{-\nu})_* \mu$ is a Haar measure and that
\[
T_{\xi} \Gamma / \xi^{\nu} T_{\xi} \Gamma \cong \Gamma / \xi^{\nu} \Gamma 
\]
has $|\xi|^{2\nu}$ elements.

It is clear that the $\R$-action on $X$ preserves the measure $\mu$. Next we need a Riemannian metric on $T \Fh\Lh$. Consider the Riemannian metric on the bundle
\[
T\C \times V_{\xi} \Gamma \times \R = \C \times \tilde{X}
\]
over $\tilde{X}$ given by:
\[
\tilde{g}_{(z , \hat{v} , t)} (\xi , \eta) = e^t \RRe (\xi \overline{\eta}) \; .
\]
It is $H$-invariant and induces a Riemannian metric $g$ on $T \Fh \Lh$ and hence on $\Lambda^{\hullet} T^* \Fh\Lh$. Together with $\mu$ we get a scalar product on $\Ah^{\hullet}_{\Fh\Lh} (X)$. With respect to this scalar product $d_{\Fh\Lh}$ has an adjoint (the formal adjoint)
\[
d^{\dagger}_{\Fh\Lh} : \Ah^{\hullet}_{\Fh\Lh} (X) \longrightarrow \Ah^{\hullet -1}_{\Fh\Lh} (X) \; .
\]
It is given by
\[
d^{\dagger}_{\Fh\Lh} = e^t d^{\dagger}_{\C} \; .
\]
Here $d^{\dagger}_{\C}$ is the formal adjoint of $d_{\C} : \Ah^{\hullet} (\C) \to \Ah^{\hullet +1} (\C)$ applied to forms in $\Ah^{\hullet +1}_{\Fh\Lh} (X)$.

Explicitly we have the following formulas:
\begin{eqnarray}
  \label{eq:7}
  d^0_{\Fh\Lh} f & = & \frac{\partial f}{\partial x} dx + \frac{\partial f}{\partial y} dy \quad \mbox{for} \; f \in \Ah^0_{\Fh\Lh} (X) = C^{\infty} (X) \\
d^{0\dagger}_{\Fh\Lh} \omega & = & -e^t \left( \frac{\partial \alpha}{\partial x} + \frac{\partial \beta}{\partial y} \right) \quad \mbox{for} \; \omega = \alpha \,dx + \beta \, dy \in \Ah^1_{\Fh\Lh} (X) \label{eq:8} \\
d^1_{\Fh\Lh} \omega & = & \left( \frac{\partial \beta}{\partial x} - \frac{\partial \alpha}{\partial y} \right) \, dx\,dy \quad \mbox{for} \; \omega = \alpha \, dx + \beta \, dy \in \Ah^1_{\Fh\Lh} (X) \label{eq:9} \\
d^{1\dagger}_{\Fh\Lh} \eta & = & e^t \left( \frac{\partial \gamma}{\partial y} \, dx - \frac{\partial \gamma}{\partial x} \, d y \right) \quad \mbox{for} \; \eta = \gamma \,dx \, dy \in \Ah^2_{\Fh\Lh} (X) \; . \label{eq:10}
\end{eqnarray}
Note that $dx$ and $dy$ do not define global forms (along $\Fh\Lh$) on $X$. The coefficients $\alpha , \beta , \gamma$ are smooth functions on $\tilde{X}$ which satisfy suitable invariance properties with respect to the $H$-action, so that $\omega , \eta$ descend to forms on $X$.

Let $\Ah^{\hullet}_{\Fh\Lh,L^2} (X)$ be the $L^2$-completion of $\Ah^{\hullet}_{\Fh\Lh} (X)$ and view $d_{\Fh\Lh}$ as an unbounded operator. We define the weak closure of $d_{\Fh\Lh}$ to be the closed unbounded operator
\[
\tilde{d}_{\FL} = d^{\dagger *}_{\FL} \; .
\]
Set
\[
\Harm^{\hullet}_{L^2} (X) = \Ker \tilde{d}_{\FL} \cap \Ker (\tilde{d}_{\FL})^* \; .
\]
This is a closed subspace of $\Ah^{\hullet}_{\FL , L^2} (X)$. For later calculations note that $(\tilde{d}_{\Fh\Lh})^* \subset \widetilde{d^{\dagger}_{\FL}}$ and therefore:
\[
\Harm^{\hullet}_{L^2} (X) \subset \widetilde{\Harm^{\hullet}_{L^2}} (X) = \Ker \tilde{d}_{\FL} \cap \Ker \widetilde{d^{\dagger}_{\FL}} \; .
\]
The metric $g$ on $T\FL$ has the property that
\[
g (T_p \phi^t (\xi) , T_p \phi^t (\eta)) = e^t g (\xi , \eta)
\]
for all $p \in X$. Since $\mu$ is $\phi^t$-invariant, it follows that for $\omega , \omega' \in \Ah^j_{\FL} (X)$ we have:
\[
(\phi^{t*} \omega , \phi^{t*} \omega') = e^{jt} (\omega , \omega') \; .
\]
Thus $\exp ( -\halb jt) \phi^{t*}$ is an orthogonal operator on the real Hilbert space $\Ah^j_{\FL} (X)$. It follows that $\phi^{t*}$ has a unique extension to $\Ah^j_{\FL, L^2} (X)$ such that $\exp (-\halb jt) \phi^{t*}$ is orthogonal. It follows that $\phi^{t*}$ commutes not only with $\tilde{d}_{\FL}$ but also with $\tilde{d}^*_{\FL}$. In particular $\phi^{t*}$ leaves $\Harm^{\hullet}_{L^2} (X)$ invariant and the representation $U_j$ of $\R$ on the real Hilbert space $\Harm^j_{L^2} (X)$ given by
\[
U_j (t) = \exp \left( -\halb jt \right) \phi^{t*}
\]
is orthogonal.

Let $\alpha \in C^{\infty}_0 (\R)$ be a test function on $\R$. Consider the bounded operator
\[
S_j (\alpha) = \int_{\R} \alpha (t) \phi^{t*} \, dt = \int_{\R} \exp \left( \halb jt \right) \alpha (t) U_j (t) \, dt
\]
on $\Harm^j_{L^2} (X)$. 

\begin{theorem}
  \label{t31}
For every $\alpha \in C^{\infty}_0 (\R)$ the operator $S_j (\alpha)$ is of trace class on $\Harm^j_{L^2} (X)$ and its trace is given for $j = 0,1,2$ by:
\begin{eqnarray*}
  \Tr (S_0 (\alpha) \tei \Harm^0_{L^2} (X)) & = & \sum_{\nu \in \Z} \Phi \left( \frac{2\pi i \nu}{\log q} \right) \\
\Tr (S_1 (\alpha) \tei \Harm^1_{L^2} (X)) & = & \sum_{\rho} \Phi (\rho) \\
\Tr (S_2 (\alpha) \tei \Harm^2_{L^2} (X)) & = & \sum_{\nu \in \Z} \Phi \left( 1 + \frac{2 \pi i \nu}{\log q} \right) \; .
\end{eqnarray*}
Here $\rho$ runs over the zeroes of $\hat{\zeta}_E (s)$. The eigenvalues of $\phi^{t*}$ on $\Harm^j_{L^2} (X)$ are the numbers
\[
\begin{array}{ll}
\dis \exp \left( \frac{2\pi i \nu t}{\log q} \right) & \mbox{for} \; \nu \in \Z \; \mbox{if} \; j = 0 \\
\exp (t\rho) & \mbox{for} \; \hat{\zeta}_E (\rho) = 0 \; \mbox{if} \; j = 1 \\
\dis \exp \left( 1+ \frac{2 \pi i \nu t}{\log q} \right) & \mbox{for} \; \nu \in \Z \; \mbox{if} \; j = 2 \; .
\end{array}
\]
\end{theorem}

Before giving the proof of the theorem we make the following observation. The maps
\[
\alpha \longmapsto \Tr (S_j (\alpha) \tei \Harm^j_{L^2} (X))
\]
define distributions $T_j$. As in the case of manifolds we define the transversal index $\Ind_t (d_{\FL})$ of the de Rham complex (5) along $\FL$ to be the alternating sum:
\[
\Ind_t (d_{\FL}) = T_0 - T_1 + T_2 \; .
\]

\begin{cor}
  \label{t32}
For the function field $F$ of an elliptic curve over $\F_q$ the left hand side of the explicit formula (2) equals $\Ind_t (d_{\FL}) (\alpha)$.
\end{cor}

\begin{proofof}
  {\bf theorem 3.1} We discuss the trace on $\Harm^1_{L^2} (X)$. The other cases are similar but easier. A $1$-form $\omega = \alpha \, dz + \beta d\overline{z}$ in $\Ah^1_{\FL,L^2} (X)_{\C}$ is given by (classes of) $V$-invariant $\C$-valued $L^2_{\loc}$-functions on $\tilde{X}$ such that for $\lambda = \nu \log q \in \Lambda$ we have
\[
\lambda^* \alpha = \xi^{\nu} \alpha \quad \mbox{and} \quad \lambda^* \beta = \xi^{\nu} \beta \; .
\]
It follows from equations (8) and (9) that $\omega$ is in $\widetilde{\Harm}^1_{L^2} (X)$ precisely when the distributional (partial) derivatives of suitable representatives $\alpha$ and $\beta$ with respect to both $x$ and $y$ vanish. The argument uses Fubini's theorem. Thus $\alpha$ and $\beta$ can be chosen to be independent of the $\C$-coordinates. 

Fourier theory shows that the action of $V$ by translation on $V_{\xi} \Gamma$ is ergodic with respect to $\mu_{\xi}$. It follows that $\alpha , \beta$ are independent of the $V_{\xi} \Gamma$-coordinate as well. 

Let $\Log_q$ denote the branch of the complex logarithm for the basis $q$ determined by $-\pi < \Arg \alpha \le \pi$. According to what we have seen, we may view:
\begin{equation}
  \label{eq:11}
  \alpha \cdot \exp (-t \Log_q \xi) \quad \mbox{and} \quad \beta \cdot \exp (-t \Log_q \oxi)
\end{equation}
as elements of $L^2 (\R / \Lambda ; \C)$.

The continuous functions being dense in $L^2 (\R / \Lambda; \C)$ it follows that $\omega$ is orthogonal to the image of $\tilde{d}_{\FL}$ i.e. that $\tilde{d}^*_{\FL} \omega = 0$. Hence we have:
\[
\Harm^1_{L^2} (X) = \widetilde{\Harm}^1_{L^2} (X) \; .
\]

It follows that $\Harm^1_{L^2} (X)_{\C}$ has an orthogonal basis consisting of the $1$-forms
\[
\exp t \left( \Log_q \xi + \frac{2 \pi i \nu}{\log q} \right)\,dz \quad \mbox{and} \quad \exp t \left( \Log_q \oxi + \frac{2 \pi i \nu}{\log q} \right)\, d\overline{z} \quad \mbox{for} \; \nu \in \Z \; .
\]
These are eigenvectors with eigenvalue $\exp (t_0 \rho)$ for the operator $\phi^{t_0*}$. Here:
\[
\rho = \Log_q \xi + \frac{2 \pi i \nu}{\log q} \quad \mbox{resp.} \quad \rho = \Log_q \oxi + \frac{2 \pi i \nu}{\log q} \quad \mbox{for} \; \nu \in \Z \; .
\]
The zeta function of $E$ is known to equal the following rational function in $q^{-s}$ c.f. \cite{S} Ch. V:
\[
\hat{\zeta}_E (s) = \frac{(1 - \pi q^{-s}) (1 - \opi q^{-s})}{(1 - q^{-s}) (1 - q^{1-s})} \; .
\]
Hence the $\rho$'s are exactly the zeroes of $\hat{\zeta}_E (s)$. The assertions on the operator $S_1 (\alpha)$ follow from straightforeward estimates.
\end{proofof}

According to corollary 3.2 the explicit formula (2) may be viewed as a formula for the transversal index $\Ind_t (d_{\Fh})$. We will now check that its right hand side can be described in terms of the closed orbits of the $\R$-action, similarly as in theorem 2.2. 

\begin{prop}
  \label{t33}
There is a natural bijection between the set of valuations $w$ of $F = \F_q (E_0)$ and the set of compact $\R$-orbits on $X$. It has the following property: If $w$ corresponds to $\gamma = \gamma_w$, then we have
\[
l (\gamma_w) = \log N (w) \; .
\]
\end{prop}

\begin{proof}
  Let $V$ act on $\C \times V_{\xi} \Gamma$ by the formula
\[
(z , \hat{v}) \cdot v = (z + v , \hat{v} - v) 
\]
and let
\[
\oM = \C \times_V V_{\xi} \Gamma
\]
be the quotient. Then $\Lambda$ acts on $\oM$ via
\[
(z , \hat{v}) \cdot \lambda = (\xi^{-\nu} z , \xi^{-\nu} \hat{v}) \quad \mbox{for} \; \lambda = \nu \log q
\]
and we may write $X$ as the suspension
\[
X = \oM \times_{\Lambda} \R \; .
\]
Identify $\oM$ with the image of $\oM \times 0$ in $X$. Then the map
\[
\gamma \longmapsto \gamma \cap \oM
\]
gives a bijection between the compact $\R$-orbits of $X$ of length $n \log q$ and the $\lambda = (- \log q)$-orbits on $\oM$ of order $n$. Note that $(-\log q) \in \Lambda$ operates by diagonal multiplication with $\xi$ on $\oM$.

The natural map
\[
\C \times_{\Gamma} T_{\xi} \Gamma \silo \C \times_V V_{\xi} \Gamma
\]
is an isomorphism and equivariant with respect to diagonal multiplication with $\xi$ on both sides.

Under the projection
\[
\C \times_{\Gamma} T_{\xi} \Gamma \longrightarrow \C / \Gamma
\]
the orbits of order $n$ of diagonal $\xi$-multiplication on $\C \times_{\Gamma} T_{\xi} \Gamma$ are mapped bijectively onto the $\xi$-orbits of order $n$ on $\C / \Gamma$. The inverse map sends the orbit of $z + \Gamma$ to the orbit of $[z , \hat{\gamma}]$ where
\[
\hat{\gamma} = (1 - \xi^n)^{-1} \gamma \quad \mbox{if} \quad \gamma = \xi^n z - z \in \Gamma \; .
\]
Note that the endomorphism $1 - \xi^n$ of $T_{\xi} \Gamma$ is invertible with inverse
\[
(1 - \xi^n)^{-1} = \sum^{\infty}_{\mu=0} \xi^{n\mu} \; .
\]
It remains to construct a bijection between the $\Theta (\xi)$-orbits on $E (\C)$ of order $n$ with the $\varphi_0$-orbits on $E_0 (\oF_q)$ of the same order. The latter are in bijection with the closed points of $E_0$ of degree $n$ i.e. with the valuations $w$ of $F = \F_q (E_0)$ of degree $n$. Note that such bijections have been constructed by Ihara in the cyclotomic, modular and Shimura curve cases \cite{I1}, \cite{I2}, \cite{I3}, \cite{I4}. However for elliptic curves there does not seem to be a reference in the literature although this case is certainly known to Ihara as well.

Let $\oL$ be the algebraic closure of $L$ in $\C$ and let $\oR$ be the integral closure of $R$ in $\oL$. Fix a maximal ideal $\oemm$ over $\emm$ in $\oR$. Then
\[
\oR / \oemm = \oF_q. 
\]
Since $\Eh$ is proper over $\spec R$ the natural map
\[
\Eh (\oR) \silo E (\oL)
\]
is an isomorphism. The torsion points of $E (\C)$ being algebraic over $L$, the inclusion
\[
E (\oL) \hookrightarrow E (\C)
\]
induces an isomorphism on the torsion subgroups. 

Consider the reduction map obtained by composition:
\[
\red : E (\C)_{\tors} \cong \Eh (\oR)_{\tors} \longrightarrow \Eh (\oR / \oemm)_{\tors} = E_0 (\oF_q) \; .
\]
It is equivariant with respect to the actions by $\Theta (\xi)$ on the left and by the $q$-Frobenius $\varphi_0$ on the right. The group
\[
(E (\C)_{\tors})^{\Theta (\xi^n)=1} = E (\C)^{\Theta (\xi^n) = 1} = (\C / \Gamma)^{\xi^n = 1} = (\xi^n - 1)^{-1} \Gamma / \Gamma
\]
has order $|\xi^n -1|^2$. The group
\[
E_0 (\oF_q)^{\varphi^n_0=1} = E_0 (\F_{q^n})
\]
has the same order, since it is known that
\[
|\xi^n-1|^2 = q^n + 1 - (\xi^n + \xi^n)= \card E_0 (\F_{q^n}) \; .
\]
For any integer $N \ge 1$ prime to $p$ the restriction
\[
\red : E (\C)_N \silo E_0 (\oF_q)_N
\]
is an isomorphism, \cite{S} IV, prop. 2.3 and VII prop. 3.1.

Now assume that $E_0$ is supersingular. Then $E_0 (\oF_q)$ has no $p$-torsion. Hence $N = |\xi^n-1|^2 = \card E_0 (\F_{q^n})$ is prime to $p$. It follows that reduction induces an isomorphism:
\[
\red : E (\C)^{\Theta (\xi^n) -1} = (E (\C)_N)^{\Theta (\xi^n) -1} \silo (E_0 (\oF_q)_N)^{\varphi^n_0 = 1} = E_0 (\F_{q^n}) \; .
\]
Therefore $\red$ also induces a bijection between the $\Theta (\xi)$-orbits on $E (\C)$ of order $n$ and the $\varphi_0$-orbits on $E_0 (\oF_q)$ of the same order.

We now assume that $E_0$ is ordinary. In this case we had taken $\Eh$ to be the canonical lifting of $E_0$ to $R = W (\F_q)$. Then the canonical sequence of $p$-divisible groups
\[
0 \longrightarrow \G_m (p) \longrightarrow \Eh (p) \longrightarrow \Q_p / \Z_p \longrightarrow 0 
\]
over $R$ is split. Hence we get a commutative diagram with compatible splittings:
\[
\xymatrix{
0 \ar[r] & \mu_{p^n} (\oR) \ar[r] \ar[d] & \Eh_{p^n} (\oR) \ar[r] \ar[d]^{\red} & \Z / p^n  \ar@/_/[l] \ar[r] \ar@{=}[d] & 0 \\
0 \ar[r] & \mu_{p^n} (\oF_q) \ar[r] & (E_{0 p^n}) (\oF_q) \ar[r] & \Z^{p^n}  \ar@/_/[l] \ar[r] & 0 .
}
\]

This in turn gives a split exact sequence:
\[
0 \longrightarrow \mu_{p^n} (\oR) \longrightarrow E (\oL)_{p^n} \xrightarrow{\red} E_0 (\oF_q)_{p^n} \longrightarrow 0
\]
and hence, by the above, a split exact sequence for any $N \ge 1$
\[
0 \longrightarrow \mu_{N_p} (\oR) \longrightarrow E (\oL)_N \xrightarrow{\red} E_0 (\oF_q)_N \longrightarrow 0 \; .
\]
Here $N_p$ is the largest power of $p$ dividing $N$.

Taking
\[
N = |\xi^n -1|^2 = \card E_0 (\F_{q^n})
\]
as above and passing to $\Theta (\xi)^n - 1$ resp. $\varphi^n_0 -1$ fixed-modules we get the exact ($\red$ was split!) sequence:
\[
0 \longrightarrow \mu_{N_p} (\oR)^{\Theta (\xi)^n -1} \longrightarrow E (\oL)^{\Theta (\xi)^n -1} \xrightarrow{\red} E_0 (\F_{q^n}) \longrightarrow 0 \; .
\]
Since the middle and right hand group have the same order it follows again that
\[
\red : E (\C)^{\Theta (\xi)^n -1} = E (\oL)^{\Theta (\xi)^n -1} \silo E_0 (\F_{q^n})
\]
is an isomorphism. We conclude as before.

The fact that $\mu_{N_p} (\oR)^{\Theta (\xi)^n - 1}$ vanishes can be seen directly.

Namely $\Theta (\xi)$ acts on $\zeta \in \mu_{p^{\infty}}$ by raising $\zeta$ to the $q$-th power. But if $\zeta^{q^n} = \zeta$ for some $n \ge 1$, if follows that $\zeta = 1$. 
\end{proof}

Combining corollary 3.2 and proposition 3.3 we arrive at the following index theoretic way to write the explicit formula (2) for an ordinary elliptic curve:

\begin{cor}
  \label{t34}
The following equality holds in $\Dh' (\R)$:
\[
\Ind_t (d_{\FL}) = \sum_{\gamma} l (\gamma) \sum_{k \ge 1} \delta_{kl (\gamma)} + \sum_{\gamma} l (\gamma) \sum_{k \le -1} e^{kl (\gamma)} \delta_{kl (\gamma)} \; .
\]
\end{cor}
This looks exactly like the transversal index formula 2.2, except for the factor $e^{kl (\gamma)}$ in front of the Dirac distribution $\delta_{kl (\gamma)}$ for $k \le -1$. Geometrically this factor results from the conformal behaviour of our metric $g$ under the flow -- a behaviour which is linked to the solenoidal structure of $X$.

In our case $\chi_{\Co} (\Fh , \mu)$ can be shown to equal $\chi (E (\C)) \log q$ which is zero.

It is possible to prove 3.4 directly using e.g. the Poisson summation formula as in the next section. Our example suggests that under suitable conditions transversal index theory generalizes to solenoids or even more general laminated spaces instead of manifolds. This is also confirmed by the example in the next section.

One cannot help wondering whether the explicit formula (1) for the Riemann zeta function might have a similar interpretation as an index theorem. The arguments in \cite{D1}, \cite{D2} suggest that in this case the $\R$-action will have a fixed point and that a generalization of transversal index theory may be required where the operator is transversally elliptic except in isolated points.

We end this section by a closer inspection of how the periodic orbits are distributed within the leaves of the $\Lh$-foliation on $X$.

For $\hat{v} \in V_{\xi} \Gamma$ one checks that
\[
\Lambda_{\hat{v}} = \{ \nu \log q \tei (\xi^{-\nu} -1) \hat{v} \in V \} \subset \Lambda
\]
is a subgroup. It acts on $\C \times \R$ by the formula:
\[
(z , t) \cdot \lambda = (\xi^{-\nu} z + (\xi^{-\nu} -1) \hat{v} , t + \lambda) \quad \mbox{for} \; \lambda = \nu \log q \in \Lambda_{\hat{v}} \; .
\]
Here we use the embeddings
\[
V_{\xi} \Gamma \hookleftarrow V \hookrightarrow \C
\]
constructed above. Set $\C \times_{\Lambda_{\hat{v}}} \R = (\C \times \R) / \Lambda_{\hat{v}}$. The map
\[
\C \times_{\Lambda_{\hat{v}}} \R \hookrightarrow X \; , \; [z,t] \longmapsto [z , \hat{v} , t]
\]
is an $\R$-equivariant diffeomorphism onto the $\Lh$-leaf obtained by projecting $\C \times \{ \hat{v} \} \times \R$ to $X$. 

For $\Lambda_{\hat{v}} = 0$ we have $\C \times_{\Lambda_v} \R = \C \times \R$. For $\Lambda_{\hat{v}} \neq 0$ we may view $\C \times_{\Lambda_{\hat{v}}} \R$ as a complex line bundle over the circle $\R / \Lambda_{\hat{v}}$ with an $\R$-action.

Note that since the set
\[
\bigcup_{\nu \in \Z \ohne 0} (\xi^{-\nu} -1)^{-1} V
\]
ist countable there are at most countably many $\hat{v}$ with $\Lambda_{\hat{v}} \neq 0$. Hence only countably many $\Lh$-leaves are diffeomorphic to $\C \times_{\Lambda_{\hat{v}}} \R$. The space of $\Lh$-leaves is $V_{\xi} \Gamma / H$.

It is clear that $\Lh$-leaves corresponding to $\hat{v} \mod H$ with $\Lambda_{\hat{v}} = 0$ do not contain compact orbits of the $\R$-action.

Assume therefore that $\Lambda_{\hat{v}} \neq 0$ and let $\lambda_0 = \nu_0 \log q$ be a generator of $\Lambda_{\hat{v}}$. Thus $\nu_0 = [\Lambda : \Lambda_{\hat{v}}]$.

The two embeddings $V_{\xi} \Gamma \hookleftarrow V \hookrightarrow \C$ being $K$-equivariant, it follows that for all $k \in \Z \ohne 0$ we have:
\[
(\xi^{-\nu_0 k} -1)^{-1} (\xi^{-\nu_0 k}-1) \hat{v} = (\xi^{-\nu_0} -1)^{-1} (\xi^{-\nu_0} -1) \hat{v} \; .
\]
It follows that the only finite orbits of the $\Lambda_{\hat{v}}$-action on $\C$ given by
\[
z \cdot \lambda = \xi^{-\nu} z + (\xi^{-\nu} -1) \hat{v} \quad \mbox{for} \; \lambda = \nu \log q \in \Lambda_{\hat{v}} 
\]
consists of the single point $z_0 = - (\xi^{-\nu_0} -1) \hat{v}$.

Hence the suspended flow, which is just the $\R$-action on $\C \times_{\Lambda_{\hat{v}}} \R$ has precisely one compact orbit $\gamma_{\hat{v}}$ and it has length
\[
l (\gamma_{\hat{v}}) = \lambda_0 = [\Lambda : \Lambda_{\hat{v}}] \log q = \vol \, \R / \Lambda_{\hat{v}} \; .
\]


\section{A holomorphic transversal index formula on $X$}
\label{sec:4}

We keep the notations from the last section. The foliation $\FL$ of $X$ carries an obvious complex structure and we can consider the Dolbeault complex along $\FL$
\[
0 \longrightarrow \Ah^{0,0}_{\FL} (X) \xrightarrow{\overline{\partial}_{\Fh}} \Ah^{0,1}_{\FL} (X) \longrightarrow 0 \; .
\]
As the hermitian metric on $T_c \FL$ etc. we take the metric induced by
\[
h_{(z , \hat{v} , t)} (\xi , \eta) = e^t \xi \overline{\eta} \; .
\]
Looking at Lazarov's result 2.3 and taking the factor $e^{k l (\gamma)}$ in 3.4 into account makes the statement of the next result plausible. 

Note also that in our case $\chi_{\Co} (\Fh , \Oh , \mu)$ can be shown to equal\\
$\chi (E (\C) , \Oh) \log q$ which is zero.

The precise definition of $\Ind_t (\overline{\partial}_{\FL})$ is given in the proof of 4.1.

\begin{theorem}
  \label{t41}
In $\Dh' (\R)$ we have:
\begin{eqnarray*}
  \Ind_t (\overline{\partial}_{\FL}) & = & \sum_{\gamma} l (\gamma) \sum_{k \ge 1} \ddet_{\C} (1 - T_{cx} \phi^{kl (\gamma)} \tei T_{cx} \FL)^{-1} \delta_{kl (\gamma)} \\
 & & + \sum_{\gamma} l (\gamma) \sum_{k \le -1} e^{kl (\gamma)} \ddet_{\C} (1 - T_{cx} \phi^{kl (\gamma)} \tei T_{cx} \FL)^{-1} \delta_{kl (\gamma)} \; .
\end{eqnarray*}
Here $\gamma$ runs over the compact $\R$-orbits on $X$ and in the $k$-sums, $x$ is any point on $\gamma$.
\end{theorem}

\begin{proof}
  For $f \in \Ah^{0,0}_{\FL} (X)_{\C} = C^{\infty} (X , \C)$ we have
  \begin{equation}
    \label{eq:15}
    \overline{\partial}_{\FL} f = \frac{\partial f}{\partial \overline{z}} d\overline{z}
  \end{equation}
and for $\omega = \beta \, d\overline{z}$ in $\Ah^{0,1}_{\FL} (X)$
\begin{equation}
  \label{eq:16}
  \opartial^{\dagger}_{\FL} \omega = -2 e^t \frac{\partial \beta}{\partial z} \; .
\end{equation}
Define $\tilde{\opartial}_{\FL}$ and $(\tilde{\opartial}_{\FL})^*$ as closed unbounded operators on $\Ah^{0,0}_{\FL , L^2} (X)_{\C}$ and $\Ah^{0,1}_{\FL , L^2} (X)_{\C}$ as before.

Let $f$ be in the kernel of $\tilde{\opartial}_{\FL}$. By elliptic regularity $f$ can be represented by a $\C$-valued $H$-invariant $L^2_{\loc}$-function on $\tilde{X}$ which we also denote by $f$. The function $f_t$ on $\oM = \C \times_V V_{\xi} \Gamma$ defined by
\[
f_t [z , \hat{v}] = [z , \hat{v} , t]
\]
is in $\Lh^2 (\oM)$ for Lebesgue almost all $t$ in $\R$. This follows from Fubini's theorem since $f$ is in $\Lh^2$. We now observe that $\oM$ is a compact abelian group with character group
\[
\oM^{\vee} = \{ \chi \otimes \chi' \tei \chi \in \C^{\vee} , \chi' \in (V_{\xi} \Gamma)^{\vee} \; \mbox{such that} \; \chi \, |_V = \chi' \, |_V \} \; .
\]
For $f_t$ in $\Lh^2 (\oM)$ we may therefore write
\begin{equation}
  \label{eq:17}
  f_t (z , \hat{v}) = \sum_{\chi \, |_V = \chi' \, |_V} a_{\chi, \chi'} \chi (z) \chi' (\hat{v}) \quad \mbox{in} \; L^2 (\oM) \; .
\end{equation}
Every character $\chi$ of $\C$ can be written in the form
\[
\chi (z) = \chi_w (z) = \exp (wz - \overline{w} \overline{z})
\]
for some unique $w \in \C$. By assumption $\opartial_z f_t = 0$ in the distributional sense. Since distributional derivatives commute with convergent series of distributions it follows that
\[
\opartial_z f_t (z , \hat{v}) = \sum_{\chi \,|_V = \chi' \, |_V} a_{\chi , \chi'} \frac{\partial \chi (z)}{\partial \overline{z}} \chi' (\hat{v}) \, d\overline{z} = 0 \; .
\]
Since
\[
\frac{\partial \chi_w}{\partial \overline{z}} = - \overline{w} \chi_w \; ,
\]
it follows that in (17) we have $a_{\chi , \chi'} = 0$ for all $\chi \neq 1$. But a character $\chi'$ of $V_{\xi} \Gamma$ with $\chi' = 1$ on $V$ is trivial since $V$ is dense in $V_{\xi} \Gamma$. It follows that $f_t$ is constant in $L^2 (\oM)$. 

Hence we have seen that we have an $\R$-equivariant isomorphism:
\[
\Ker \tilde{\opartial}_{\FL} = L^2 (\R / \Lambda , \C) \; .
\]
It follows that for every $\alpha \in C^{\infty}_0 (\R)$ the operator
\[
S_0 (\alpha) = \int_{\R} \alpha (t) \phi^{t*} \, dt \quad \mbox{on} \; \Ker \tilde{\opartial}_{\FL}
\]
is of trace class. Its trace is given by the formula
\[
T_0 (\alpha) = \Tr (S_0 (\alpha) \tei \Ker \tilde{\opartial}_{\FL}) = \sum_{\nu \in \Z} \Phi \left( \frac{2\pi i\nu}{\log q} \right) \; .
\]
Similar arguments as before show that $\Ker \widetilde{\opartial^{\dagger}} \supset \Ker \tilde{\opartial}^*_{\FL}$ may be identified with
\[
L^2 (\R / \Lambda , \C) \cdot \exp (t \Log_q \oxi) \, d\overline{z}
\]
which in fact lies in the kernel of $\tilde{\opartial}^*$. Hence we have
\[
\Ker \widetilde{\opartial^{\dagger}} = \ker \tilde{\opartial}^* 
\]
in our case. Moreover the operator
\[
S_1 (\alpha) = \int_{\R} \alpha (t) \phi^{t*} \, dt \quad \mbox{on} \; \Ker \tilde{\opartial}^*_{\FL}
\]
is of trace class and
\[
T_1 (\alpha) = \Tr (S_1 (\alpha) \tei \Ker \tilde{\opartial}^*_{\FL}) = \sum_{\nu \in \Z} \Phi \left( \Log_q \oxi + \frac{2 \pi i \nu}{\log q} \right) \; .
\]
Both $T_0$ and $T_1$ are distributions on $\R$ and we set
\[
\Ind_t (\opartial_{\FL}) = T_0 - T_1 \; .
\]
For $\psi \in C^{\infty}_0 (\R)$ and
\[
\hat{\psi} (x) = \int_{\R} e^{-2 \pi ixt} \psi (t) \, dt
\]
the Poisson summation formula asserts that
\[
\sum_{n \in \Z} \psi (n) = \sum_{n \in \Z} \hat{\psi} (n) \; .
\]
It implies that
\[
\Ind_t (\opartial_{\FL}) = \log q \sum_{n \in \Z} (1 - \oxi^n) \delta_{n \log q} \; .
\]
Using the formulas valid for $n \ge 1$:
\[
\card E_0 (\F_{q^n}) = ( 1- \xi^n) (1 - \oxi^n)
\]
and 
\[
\card E_0 (\F_{q^n}) = \sum_{\deg w \tei n} \deg w
\]
we find after a short calculation:
\begin{eqnarray*}
  \Ind_t (\opartial_{\FL}) & = & \sum_w \log Nw \sum_{k \ge 1} (1 - \xi^{k \deg w})^{-1} \delta_{k \log N (w)} \\
& & + \sum_w \log N w \sum_{k \le -1} q^{k \deg w} (1 - \xi^{k \deg w})^{-1} \delta_{k \log N (w)} \; .
\end{eqnarray*}
Now if $x$ is a point on a closed orbit $\gamma$ then the map induced by $\phi^{kl (\gamma)}$ on the complex tangent space $T_{cx} \FL$ is multiplication by $\xi^{k (\log q)^{-1} l (\gamma)}$.

Using proposition 3.3 we therefore get the assertion.
\end{proof}

\noindent
Mathematisches Institut\\
Westf. Wilhelms-Universit\"at\\
Einsteinstr. 62\\
48149 M\"unster\\
Germany\\
deninge@math.uni-muenster.de
\end{document}